\documentclass[10pt]{amsart}
\usepackage{amssymb,bbold}
\usepackage[all]{xy}

\theoremstyle{plain}
\newtheorem{theorem}{Theorem}

\theoremstyle{definition}

\title[Continuous functions that cut the real axis very often]{Continuous functions that cut the real axis very often}
\author{Omid Zabeti}
\address{O. Zabeti, Department of Mathematics, Faculty of Mathematics, University of Sistan and Baluchestan,  P.O. Box 98135-674, Zahedan (Iran)}
\email{o.zabeti@gmail.com}
\begin{document}
\begin{center}
\bf{Dedicated to Miguel Lacruz, my friend and my colleague}
\end{center}
\begin{abstract}
We consider  continuous functions \(f \colon [0,1] \to {\mathbb R}\) that cut the real axis at every point of a measurable set  of positive measure and we construct  examples where \(f\) fails to have bounded variation,  and at the opposite end, where \(f\) admits derivatives of all orders.
\end{abstract}
\date{\today}
\maketitle
We say that  a function  \(f \colon [0,1] \to \mathbb R\) {\em cuts the real axis at a point} \(x \in [0,1]\) provided that \(f(x)=0,\) and for every \(\delta >0\) there exist \(x_1,x_2 \in [0,1]\) such that \(|x-x_i| <\delta\)  for \(1 \leq i \leq 2\) and \(f(x_1) < 0 < f(x_2).\)

Recall that a function \(f \colon [0,1] \to \mathbb R\) is said to have {\em bounded variation} provided that
\[
V_0^1 (f) \colon = \sup_P  \sum_{i=1}^n |f(t_i)-f(t_{i-1})| < \infty,
\]
where the supremum is taken over all  partitions  \(P=\{0=t_0 < t_1 < \cdots <t_n=1\}\) of the interval \([0,1].\) The above expression is called the {\em total variation of \(f.\)} It is a standard fact that any differentiable function \(f\) with continuous derivative \(f^\prime\) has bounded variation and moreover,
\[
V_0^1(f)= \int_0^1 |f^\prime (x)|\, dx.
\]

The purpose of this paper is to construct  an  example  of a continuous function that cuts the real axis at every point of a measurable set  of positive measure and that  fails to have bounded variation, and at the opposite end, another example of a function  that cuts the real axis at every point of a measurable set  of positive measure and that admits derivatives of all orders.
				
Our first construction is based on the properties of the generalized Cantor sets.  We follow the discussion of the generalized Cantor sets in the book of Folland \cite[p.40]{folland}. Let \((\xi_n)\) be a sequence of positive numbers such that \(\xi_0=1\) and \(\xi_n > 2 \xi_{n+1}.\) Remove from \([0,1]\) the open middle interval of length \(\xi_0 - 2 \xi_1,\) obtaining a set \(B_1\) that is the union of two disjoint closed intervals \([0, \xi_1]\) and \([1-\xi_1,1].\) Proceeding inductively, having constructed \(B_n,\) remove  from each of its  \(2^n\) constituent intervals of length \(\xi_n\) the open middle interval of length \(\xi_n -2 \xi_{n+1}\) to obtain a set \(B_{n+1}\) that is the union of \(2^{n+1}\) disjoint closed intervals of length \(\xi_{n+1}.\)  Finally, the intersection
\[
B= \bigcap_{n \in \mathbb N} B_n
\]
is called a {\em generalized Cantor set.}  It is clear that \(|B|= \lim 2^n \xi_n.\) In particular, if \(0 \leq \alpha < 1\) and we take \(\xi_n = \alpha 2^{-n} + (1-\alpha) 3^{-n}\) then it is easy to check that the condition \(\xi_n > 2\xi_{n+1}\) is satisfied and \(|B|=\alpha.\) From now on we restrict our attention to this particular choice.

\begin{theorem} Let \(0 < \alpha < 1\) and let \(B\) be a generalized Cantor set with \(|B|=\alpha.\) There exists a continuous function  \(f \colon [0,1] \to \mathbb R\) that fails to have bounded variation and that cuts the real axis at every \(x \in B.\)
\end{theorem}
\begin{proof} We construct a sequence of functions \((f_n)\) that converges uniformly on \([0,1]\) to the desired function.  We proceed by induction. Let \((c_n)\) be a sequence of positive numbers such that
\[
\sum_{n=1}^\infty c_n < \infty, \quad \text{and} \quad \lim_{n \to \infty} 2^n c_n = \infty
\]
(for instance, the sequence \(c_n=n/2^n\) does the job). First of all, define  a   function  \(h_1 \colon [0, 1- 2\xi_1] \to \mathbb R\)   by the expression \(h_1(x) = c_1\sin(2\pi x/(1-2\xi_1)).\) Next, define a  function \(f_1 \colon [0,1] \to \mathbb R\) by the expression
\[
f_1(x) = 	\left \{\begin{array}{rl} 0, & \text{ if } 0 \leq x \leq \xi_1,\\ h_1(x-\xi_1), & \text{ if } \xi_1 \leq x \leq 1-\xi_1,\\ 0, & \text{ if } 1-\xi_1 \leq x \leq 1. \end{array}	
		\right .
\]
Suppose that we have constructed \(f_n\)  and  define   a function \(h_{n+1} \colon [0, \xi_n - 2\xi_{n+1}] \to \mathbb R\)  by the expression
\[
h_{n+1}(x) =c_{n+1}\sin \left ( \frac{2\pi x}{\xi_n-2\xi_{n+1}} \right ).
\]
Then, define a  continuous function \(f_{n+1} \colon [0,1] \to \mathbb R\) as follows. First, set \(f_{n+1}(x)=f_n(x)\) for all \(x \in [0,1] \backslash B_n.\)  Next,  define \(f_{n+1}\) on the set \(B_n\)  by setting its value to zero on each constituent interval of \(B_{n+1},\) and by making \(2^n\) copies of \(h_{n+1}\) on each of the open middle intervals removed from  the  constituent intervals of \(B_n.\) It follows from the construction that each \(f_n\) is a continuous function that changes sign on each   interval of \([0,1] \backslash B_n\) and that vanishes on \(\partial B_n.\) Moreover, it is clear that \(\| f_{n+1}-f_n \|_\infty  \leq  c_n+c_{n+1}\) for every \(n \geq 1,\)  so that
\[
\sum_{n=1}^\infty \|f_{n+1}-f_n\|_\infty < \infty.
\]
Hence,  the sequence \((f_n)\) converges uniformly to a continuous function, say \(f.\) We claim that \(f\) fails to have bounded variation. Indeed,  if we consider the variation of \(f\) on each of the \(2^n\) open middle intervals removed from  the  constituent intervals of \(B_n\) then we get a  lower bound on  the total variation of \(f,\)  namely, \(V_0^1(f) \geq 2^n V_0^{\xi_n-2\xi_{n+1}}(f_{n+1}) = 2^{n+2} c_{n+1} \to \infty.\) Also, we claim that \(f\) cuts the real axis at every \(x \in B \cap(0,1).\) Indeed,  since \(f\) vanishes at the endpoints of each constituent interval of  \(B_n,\)    since these endpoints are dense in \(B\) and since  \(f\) is continuous, it is clear that   \(f(x)=0.\) Now, let \(\delta >0\) and choose \(n \in \mathbb N\) large enough, so that \((x-\delta,x+\delta)\) contains one of the   open middle intervals removed from  the  constituent intervals of \(B_n.\) Since \(f_n\) changes sign on that interval, it follows that  there exist \(x_1,x_2 \in [0,1]\) such that \(|x-x_1|< \delta\) and such that    \(f(x_1)=f_n(x_1)<0<f_n(x_2)=f(x_2).\)
\end{proof}

\begin{theorem} There exists a function  \(f \colon [0,1] \to \mathbb R\) that  admits derivatives of all orders, and there exists a measurable set \(E \subseteq [0,1]\) such that \(|E|>0\) and  \(f\) cuts the real axis at every \(x \in E.\)
\end{theorem}
\begin{proof} Let \((r_n)\) denote the sequence of all rational numbers in \((0,1)\), let \(0 <\varepsilon <1,\) and for every \(n \in \mathbb N\) consider the open interval
\[
I_n= \left (r_n - \frac{\varepsilon}{2^{n+1}}, r_n +\frac{\varepsilon}{2^{n+1}} \right ).
\]
Now, consider the dense open set
\[
G = \bigcup_{n \in \mathbb N} I_n.
\]
It follows from the subadditivity of the Lebesgue measure that \(|G| \leq \varepsilon.\) Then, consider the closed set \(F= [0,1] \backslash G\) and notice that \(|F| \geq 1-\varepsilon.\)
Now,  \(G\) is a countable union of disjoint open intervals, say
\[
G= \bigcup_{n \in \mathbb N} J_n,
\]
where \(J_n=(a_n - \varepsilon_n,a_n+\varepsilon_n).\) Notice that
\[
2 \sum_{n=1}^\infty \varepsilon_n = |G| \leq \varepsilon
\]
so that in particular \(\lim \varepsilon_n=0.\) Next, let \(E \subseteq F\) be the set obtained  after removing  from \(F\)  the endpoints of  the intervals \(J_n,\) so that we still have \(|E| \geq 1-\varepsilon.\)
It is a standard fact that  the  function  \(h \colon \mathbb R \to \mathbb R\) defined by the expression
\[
h(x) = \left \{
		\begin{array}{rl}
		\exp \left (\displaystyle{ -\frac{1}{1-x^2}} \right ), & \text{ if } |x|<1,\\
		0, & \text{ otherwise}
		\end{array}
	\right .
\]
admits derivatives of all orders  supported on the interval \([«-1,1].\)
Then, consider the sequence of functions \(f_n \colon [0,1] \to \mathbb R\) defined by the expression
\[
f_n(x) = c_n h \left ( \frac{x-a_n}{\varepsilon_n} \right )  \sin \pi  \left ( \frac{x-a_n}{\varepsilon_n} \right )
\]
where \((c_n)\) is a sequence of positive numbers to be chosen later on.
Notice that each \(f_n\) is supported on the closure of \(J_n\) and  its graph cuts the real axis at \(a_n.\) It follows from Leibniz's rule for the derivatives of higher order that
\[
f^{(p)}_n (x) = c_n \varepsilon_n^{-p} \pi^{p}  \sum_{k=0}^p \binom{p}{k}h^{(k)}  \left ( \frac{x-a_n}{\varepsilon_n} \right ) \pi^{-k}  \sin^{(p-k)} \pi  \left ( \frac{x-a_n}{\varepsilon_n} \right ).
\]
Notice that  \( \|f_n^{(p)}\|_\infty \leq c_n \varepsilon_n^{-p} M_p,\) where
\[
M_p = \pi^p \sum_{k=0}^p \binom{p}{k}  \pi^{-k} \|h^{(k)}\|_\infty.
\]

\noindent
Next, consider the function \(f \colon [0,1] \to \mathbb R\)  defined by the infinite series of functions
\[
f(x)=\sum_{n=1}^\infty f_n (x),
\]
and  choose a sequence \((c_n)\) such that, for every \(p \in \mathbb N,\)
\[
\sum_{n=1}^\infty c_n \varepsilon_n^{-p} < \infty.
\]
It is an easy exercise to show that the sequence \(c_n=\exp(-1/\varepsilon_n)\) satisfies the above condition.
A  standard application of the Weierstrass M-test leads  to the conclusion that \(f\) admits derivatives of all orders  that are represented by the uniformly convergent series
\[
f^{(p)} (x)=\sum_{n=1}^\infty f_n^{(p)} (x).
\]
Finally, we must show that  \(f\) cuts the real axis at every \(x \in E.\) It is clear that \(f(x)=0,\) since each \(f_n\) is supported on the closure of the interval  \(J_n\) and \(x\) is not in the closure of any of these intervals. Since \(G\) is dense in \([0,1],\)  there exists a sequence \((x_j)\) in \(G\) that converges to \(x.\) Now, for every \(j \in \mathbb N\) there exists \(n_j \in \mathbb N\) with \( |x_j - a_{n_j}| < \varepsilon_{n_j}. \)  Then, we claim that the set of indices \(\{n_j \colon j \in \mathbb N\}\) is infinite, for otherwise one of the indices is repeated infinitely often, that is, there exists \(j_0 \in \mathbb N\) and there exists an increasing sequence \((j_k)\) such that   \(n_{j_k}=n_{j_0}\) for all \(k \in \mathbb N.\) Thus, the subsequence \((x_{j_k})\) stays in the interval \(J_{n_{j_0}},\) and since \(x=\lim x_{n_{j_k}}\) and \(x  \notin J_{n_{j_0}},\) it follows that \(x\) is an endpoint of  \(J_{n_{j_0}}.\) A contradiction has arrived, for \(E\) is the set obtained from \(F\) after removing the endpoints of the intervals \(J_n.\) Since the set of indices \(\{n_j \colon j \in \mathbb N\}\) is infinite, there exists an increasing sequence \((j_k)\) such that the sequence \((n_{j_k})\) is also increasing. Finally, since  \(\lim \varepsilon_n=0,\) it follows that \  \((a_{n_{j_k}})\) converges to \(x.\)  Let \(\delta>0\) and choose \(k_1 \in \mathbb N\)  such  that \(a_{n_{j_k}}  \in (x-\delta,x+\delta)\) for all \(k \geq k_1.\) Then, choose \(k_2 \in \mathbb N\) such that \(\varepsilon_{n_{j_k}} < \delta \) for all \(k \geq k_2.\) Now, set \(k_0 = \max \{k_1,k_2\}\)  and notice that  \(J_{n_{j_{k_0}} } \subseteq  (x-\delta, x+ \delta ).\)
Since the function \(f_{n_{j_{k_0}}}\) cuts the real axis at \(a_{n_{j_{k_0}}},\) there are \(x_1,x_2 \in J_{n_{j_{k_0}}}\) such that
\[
f(x_1) =f_{n_{j_{k_0}}}(x_1)<0<f_{n_{j_{k_0}}}(x_2)=f(x_2),
\]
as we wanted.
\end{proof}
{\bf{acknowlegement}}. This note would not exist without inspiring and invaluable suggestions of Professor Miguel Lacruz, my friend and my colleague.

\end{document}